\documentclass{l4dc2021} 

\title[Safe Learning from Short Trajectories]{Safely Learning Dynamical Systems from Short Trajectories}
\usepackage{times}
\usepackage{mathtools}

\DeclareGraphicsExtensions{.pdf,.eps}

\ifdefined\ISSUBMISSION
\newcommand{\myproof}[1]{}
\newcommand{\myproofof}[2]{}

\else
\newcommand{\myproof}[1]{\begin{proof}#1\end{proof}}
\newcommand{\myproofof}[2]{\begin{proof}[#1]#2\end{proof}}

\fi

\author{\Name{Amir Ali Ahmadi} \Email{aaa@princeton.edu}\\
 \Name{Abraar Chaudhry} \Email{azc@princeton.edu}\\
 \addr ORFE, Princeton University
 \AND
 \Name{Vikas Sindhwani} \Email{sindhwani@google.com}\\
 \Name{Stephen Tu} \Email{stephentu@google.com}\\
 \addr Robotics at Google, New York}

\makeatletter
\ifjmlrutilssubfloats
\@ifl@t@r\fmtversion{2019/08/22}%
{
\renewcommand*\@subfigurelabel[3]{#1\subfigurelabel{#2}}
\renewcommand*\@subfigref[1]{%
{%
\def\@subfigurelabel##1##2##3{\subfigurelabel{##2}}%
\ref{#1}%
}%
}
}%
{
\renewcommand*\@subfigurelabel[2]{#1\subfigurelabel{#2}}
\renewcommand*\@subfigref[1]{%
{%
\def\@subfigurelabel##1##2{\subfigurelabel{##2}}%
\ref{#1}%
}%
}
}
\makeatother

\newcommand{\R}{\mathbb{R}}
\newcommand{\Tr}{\text{Tr}}
\newcommand{\spanof}[1]{\textnormal{span}(#1)}

\newcommand{\norm}[1]{\| #1 \|}

\newcommand{\trueA}{A_{\star}}
\newcommand{\truef}{f_{\star}}
\newcommand{\trueg}{g_{\star}}
\newcommand{\defn}{\mathrel{\mathop :}=}

\begin{document}

\maketitle

\begin{abstract}%
A fundamental challenge in learning to control an unknown dynamical system
is to reduce model uncertainty by making measurements while maintaining safety. 
In this work, we formulate a mathematical definition of what
it means to safely learn a dynamical system 
by sequentially deciding where to initialize the next trajectory.
In our framework, the state of the system is required to stay within a given safety
region under the (possibly repeated) action of all dynamical systems that are consistent with
the information gathered so far.
For our first two results, we consider the setting of safely learning linear dynamics.
We present a linear programming-based
algorithm that either safely recovers the true dynamics from trajectories of length
one, or certifies that safe learning is impossible.
We also give an efficient semidefinite representation of the
set of initial conditions whose resulting trajectories of length two
are guaranteed to stay in the safety region.
For our final result, we study the problem of safely learning a nonlinear
dynamical system. We give a second-order cone programming
based representation of the
set of initial conditions that are guaranteed to remain in the safety region
after one application of the system dynamics.\\


%
%

\end{abstract}


\section{Introduction and Problem Formulation}
The core task in model-based reinforcement learning \citep{yang2020data,nagabandi2018neural,singh2019learning,lowrey2018plan,venkatraman2016improved,kaiser2019model} is to estimate---from a small set of sampled trajectories---an unknown dynamical system prescribing the evolution of an agent's state given the current state and control input.
During the initial stages of learning, deploying even a conservative feedback policy on a real robot is fraught with risk, even if the policy achieves high task performance and safe behavior in simulation. How should the robot be ``set loose" in the real world so that the dynamics may be precisely estimated by observing state transitions, but with strong guarantees that the robot will remain safe? This interplay between \emph{safety and uncertainty while learning dynamical systems} is the central theme of this paper. 

We view the agent armed with a fixed feedback policy in closed loop over a short duration as an unknown discrete-time dynamical system
\begin{equation}\label{eq:dynamics}
x_{t+1} = \truef(x_t).
\end{equation}
We consider the problem of safe data acquisition for estimating the unknown map $\truef: \R^n \rightarrow \R^n$ from a collection of length-$T$ trajectories $\{\phi_{\truef, T}(x_k)\}_{k=1}^{m}$,
where 
$\phi_{f,T}(x) \defn ( x, f(x), \dots, f^{(T)}(x) )$. In our setting, we are given as input a set $S \subset \R^n$, called the \emph{safety region}, in which the state should remain throughout the learning process.  We say that a state $x\in\R^n$ is \emph{$T$-step safe} under a map $f: \R^n \rightarrow \R^n$ if $x$ belongs to the set $$S^T(f) \defn \{ x\in S \mid f^{(i)}(x) \in S, i=1, \ldots, T\}.$$
In order to safely learn $\truef$, we require that measurements are made only at points $x\in\R^n$ for which $x\in S^T(\truef)$. 
Obviously, if we make no assumptions about $\truef$, this task is impossible. We assume, therefore, that the map $\truef$ belongs to a set of dynamics $U_0$, which we call the \emph{initial uncertainty set}. As experience is gathered, the uncertainty over $\truef$ decreases. 
Let us denote the uncertainty set after the agent has observed $k$ trajectories $\{ \phi_{\truef, T}(x_j)\}_{j=1}^{k}$ by,
$$U_k \defn \{f\in U_0 \mid \phi_{f,T}(x_j) = \phi_{\truef,T}(x_j) \:, j=1,\dots,k \}.$$ 
For a nonnegative integer $k$, define
\[S^T_k \defn \bigcap_{f \in U_k} S^T(f) \:, \]
the set of points that are $T$-step safe under all dynamics consistent with the data after observing $k$ trajectories.
Fix a distance metric $d(\cdot, \cdot)$ over $U_0$ and a scalar $\varepsilon>0$.
Given a safety region $S\subset\R^n$ and an initial uncertainty set $U_0$, we say that \emph{$T$-step safe learning} is possible (with respect to the metric $d(\cdot, \cdot)$ and up to accuracy $\varepsilon$) if for some nonnegative integer $m$, we can sequentially choose vectors $x_1,\dots,x_m$ such that,
\begin{enumerate}
    \item \textbf{(Safety)} for each $k=1, \dots, m$, $x_{k} \in S^T_{k-1}$,
    \item \textbf{(Learning)} $\sup_{f \in U_m} d(f, \truef) \leq \varepsilon$.
\end{enumerate}
Note that for any $T'>T$, we have $S^{T'}_{k} \subseteq S^T_{k}$ for all $k$. Hence, if $T$-step safe learning is impossible, then $T'$-step safe learning is also impossible. Therefore, the highest rate of safe information assimilation during the learning process is achieved when $T=1$. One of the main contributions of this paper is to present an efficient algorithm
for the exact one-step safe learning problem (i.e., when $\varepsilon = 0$ and $T=1$) in the case where 
the dynamics in \eqref{eq:dynamics} are linear, $U_0$ is a polyhedron
in the space of $n \times n$ matrices that define the dynamics,
and $S$ is a polyhedron (\algorithmref{alg:one-step}
and \theoremref{thm:one-step}).

Suppose furthermore that initializing the unknown system at a state $x\in S$ comes at a cost of $c(x)$. 
In such a setting, we are also interested in 
safely learning at minimum measurement cost.
To do this, one naturally wants to solve an optimization problem of the type
\begin{equation}
\min_{x \in S^T_{k-1}} c(x) \:, \label{eq:greedy_safe}
\end{equation}
whose optimal solution gives us the next cheapest $T$-step safe query point $x_{k}$.
Another contribution of this paper is to 
derive exact reformulations of problem \eqref{eq:greedy_safe}, when $T \in \{1, 2\}$, in terms of tractable conic optimization problems. More specifically, under natural assumptions on $S$ and $U_0$,
when the unknown dynamics are linear,
we show that problem \eqref{eq:greedy_safe} can be formulated as a linear program when $T=1$ (\theoremref{thm:one-step LP}) and as a semidefinite program when $T=2$ (\theoremref{thm:two-step}).
Furthermore, when $T=1$ and the unknown dynamics are nonlinear (but bounded in a certain sense), we show that \eqref{eq:greedy_safe} can be formulated as
a second-order cone program (\theoremref{thm:nonlinear socp}). 
Finally, we note that we are currently preparing a draft to handle the case when $T=\infty$
using the set invariance tools of \cite{ahmadi18rdo}.

\section{Related Work}

Most related to our work is \cite{dean19safelqr},
which uses the system-level synthesis framework~\citep{anderson19sls}
to derive inner approximations to 
the infinite-step safety region of a linear system subject to
polytopic uncertainty in the dynamics and bounded disturbances.
\cite{lu17safeexploration} considers a probabilistic version of one-step safety for linear systems and also presents
an algorithm to conservatively compute the $T$-step safety regions.
Unlike these papers that focus on inner approximations of safety regions, we are able to exactly characterize one-step and two-step safety regions under our proposed framework.
We also note that we do not require any stability assumptions on
the dynamical systems we want to learn.

We also review other works focused on the general problem of safely
learning dynamics in both the control theory and reinforcement
learning literature.
\cite{berkenkamp17safeRL} combines Lyapunov functions and Gaussian process models to show how to safely explore an uncertain system and
expand an inner estimate of the region of attraction of one of its equilibrium points.
\cite{akametalu14safeGP}
uses reachability analysis to compute maximal safe regions for uncertain dynamics, and proposes Gaussian processes to iteratively refine the uncertainty. 
\cite{koller19learning}
shows how to propagate ellipsoidal uncertainty multiple steps into the future,
and utilizes this uncertainty propagation in a model predictive control framework for safely
learning to control.
\cite{wabersich18mpc}
shows how to minimally perturb a controller designed to learn a linear system in order for the system to stay within a set of constraints that guarantee reachability to a safe target set.

We also note that our work has some conceptual connections to the literature on experiment design \citep[see e.g.,][]{pukelsheim06doe, decastro19experimentdesign}.
However, this literature typically does not consider dynamical systems or notions of safety.

\section{One-Step Safe Learning of Linear Systems}\label{sec:one-step}

In this section, we focus on characterizing one-step safe learning
for linear systems.
Here, the state evolves according to 
\begin{equation}\label{eq:linear dynamics}
x_{t+1} = \trueA x_t,
\end{equation}
where $\trueA$ is an unknown $n \times n$ matrix.
We assume we know that $\trueA$ belongs to a set $U_0 \subset \R^{n \times n}$ that represents our prior knowledge of $\trueA$.
In this section, we take $U_0$ to be a polyhedron; i.e., \begin{equation}\label{eq:U_A polyhedron}
U_0 = \left\{ A \in \R^{n \times n} \mid  \Tr(V_j^T A) \leq v_j \quad j = 1, \dots ,s \right\}
\end{equation}
for some matrices $V_1,\dots,V_s \in \R^{n \times n}$ and scalars $v_1,\dots,v_s \in \R$. We also work with a polyhedral representation of the safety region $S$; i.e.,
\begin{equation}\label{eq:S polyhedron}
S = \left \{ x \in \R^n \mid  h_i^T x \leq b_i \quad i = 1, \dots ,r \right \}
\end{equation}
for some vectors $h_1,\dots,h_r \in \R^n$ and some scalars $b_1,\dots,b_r \in \R$.
We assume that making a query at a point $x\in \R^n$ comes at a cost $c^T x$, where the vector $c \in \mathbb{R}^n$ is given\footnote{In practice, measurement costs are typically nonnegative. 
If $S$ is compact for example, one can always add a constant term to 
$c^T x$ to ensure this requirement without changing any of our optimization problems.}.
An extension to more general semidefinite representable cost functions
is possible using tools of conic optimization. 

We start by finding 
the minimum cost point that is one-step safe under all valid dynamics, i.e.,\ a point $x\in S$ such that $Ax\in S$ for all $A\in U_0$. Once this is done, we gain further information by observing the action $y = \trueA x$ of system \eqref{eq:linear dynamics} on our point $x$, which further constrains the uncertainty set $U_0$. We then repeat this procedure with the updated uncertainty set to find the next minimum cost one-step safe point. 
More generally, after collecting $k$ measurements, 
our uncertainty in the dynamics reduces to the set
\begin{align}
    U_k = \{ A \in U_0 \mid  Ax_j = y_j \quad j = 1,\dots,k\}. \label{eq:uk_linear}
\end{align}
Hence, the problem of finding the next cheapest one-step
safe query point becomes:
\begin{equation}\label{eq:one-step}
\begin{aligned}
\min_{x\in\R^n} \quad & c^T x\\
\textrm{s.t.} \quad & x \in S \\
& A x \in S \quad \forall A \in U_k. 
\end{aligned}
\end{equation}


In \sectionref{sec:one-step:duality}, we show that
problem \eqref{eq:one-step} can be efficiently solved.
We then use \eqref{eq:one-step}
as a subroutine in a one-step safe learning algorithm which 
we present in \sectionref{sec:one-step:algorithm}.

\subsection{Reformulation via Duality}
\label{sec:one-step:duality}

In this subsection, we reformulate \eqref{eq:one-step} as a linear program.
To do this we introduce auxiliary variables $\mu_j^{(i)} \in \R$ and $\eta_k^{(i)} \in \R^n$ for $i=1,\dots,r$, $j=1,\dots,s$, and $k=1,\dots,m$.
\begin{theorem}\label{thm:one-step LP}
The feasible set of problem \eqref{eq:one-step} is the projection onto $x$-space of the feasible set of the following linear program:
\begin{equation}\label{eq:one-step LP}
\begin{aligned}
\min_{x,\mu,\eta} \quad & c^T x\\
\textrm{s.t.} \quad & h_i^T x \leq b_i \quad i = 1, \dots ,r \\
& \sum_{k=1}^{m} y_k^T \eta_k^{(i)} + \sum_{j=1}^{s} \mu_j^{(i)} v_j \leq b_i \quad i = 1, \dots ,r \\
& x h_i^T = \sum_{k=1}^{m} x_k \eta_k^{(i)T} + \sum_{j=1}^{s} \mu_j^{(i)} V_j^T \quad i = 1, \dots ,r \\
& \mu^{(i)} \geq 0 \quad i = 1, \dots ,r.
\end{aligned}
\end{equation}
In particular, the optimal values of \eqref{eq:one-step} and \eqref{eq:one-step LP} are the same and the optimal solutions of \eqref{eq:one-step} are the optimal solutions of \eqref{eq:one-step LP} projected onto $x$-space.
\end{theorem}
\myproof{
Using the definitions of $S$ and $U_0$, let us first rewrite \eqref{eq:one-step} as bilevel program:
\begin{equation}\label{eq:one-step bilevel}
\begin{aligned}
\min_{x} \quad & c^T x\\
\textrm{s.t.} \quad & h_i^T x \leq b_i \quad i = 1, \dots ,r \\
& \begin{bmatrix} \max_A \quad & h_i^T A x  \\ 
\textrm{s.t.} \quad & \Tr(V_j^T A) \leq v_j \quad j = 1, \dots ,s \\
& Ax_k = y_k \quad k=1, \dots, m \end{bmatrix} \leq b_i \quad i = 1, \dots ,r.
\end{aligned}
\end{equation}
We proceed by taking the dual of the $r$ inner programs, treating the $x$ variable as fixed.
By introducing dual variables $\mu_j^{(i)}$ and $\eta_k^{(i)}$ for $i=1,\dots,r$, $j=1,\dots,s$, and $k=1,\dots,m$, and by invoking strong duality of linear programming, we have
\begin{equation}\label{eq:one-step duality}
\begin{bmatrix} \max_A \: & h_i^T A x  \\ 
\textrm{s.t.} \: & \Tr(V_j^T A) \leq v_j \: j = 1, \dots ,s \\
& Ax_k = y_k \: k=1, \dots, m \end{bmatrix} = 
\begin{bmatrix} \min_{\mu^{(i)},\eta^{(i)}} \: & \sum_{k=1}^{m} y_k^T \eta_k^{(i)} + \sum_{j=1}^{s} \mu_j^{(i)} v_j \\ 
\textrm{s.t.} \: & x h_i^T = \sum_{k=1}^{m} x_k \eta_k^{(i)T} + \sum_{j=1}^{s} \mu_j^{(i)} V_j^T \\
& \mu^{(i)} \geq 0 \end{bmatrix} 
\end{equation}
for $i=1,\dots,r$.
Thus by replacing the inner problem of \eqref{eq:one-step bilevel} with the right-hand side of \eqref{eq:one-step duality}, the min-max problem \eqref{eq:one-step bilevel} becomes a min-min problem.
This min-min problem can be combined into a single minimization problem and be written as problem \eqref{eq:one-step LP}.
Indeed, if $x$ is feasible to \eqref{eq:one-step bilevel}, for that fixed $x$ and for each $i$, there exist values of $\mu^{(i)}$ and $\eta^{(i)}$ that attain the optimal value for \eqref{eq:one-step duality} and therefore the triple $(x,\mu,\eta)$ will be feasible to \eqref{eq:one-step LP}.
Conversely, if some $(x,\mu,\eta)$ is feasible to \eqref{eq:one-step LP}, it follows that $x$ is feasible to \eqref{eq:one-step bilevel}. This is because for any fixed $x$ and for each $i$, the optimal value of the left-hand side of \eqref{eq:one-step duality} is bounded from above by the objective value of the right-hand side evaluated at any feasible $\mu^{(i)}$ and $\eta^{(i)}$.
}

We remark that \eqref{eq:one-step LP}
can be modified so that
one-step safety is achieved in the presence of 
disturbances. We can
ensure, e.g.,\ using linear programming, that $A x + w \in S$
for all $A \in U_m$
and all $w$ such that $\norm{w} \leq W$, where $\norm{\cdot}$ is any
norm whose unit ball is a polytope and $W$ is a given scalar.

\subsection{An Algorithm for One-Step Safe Learning}
\label{sec:one-step:algorithm}

We start by giving a mathematical definition of (exact) safe learning specialized to the case of one-step safety and linear dynamics.
Recall the definition of the set $U_k$ in \eqref{eq:uk_linear}.
\begin{definition}[One-Step Safe Learning]
\label{def:one-step}
We say that one-step safe learning is possible if for some nonnegative integer $m$, we can sequentially choose vectors $x_k \in S$, for $k=1, \ldots, m$, and observe measurements $y_k = \trueA x_k$ such that:
\begin{enumerate}
    \item \textbf{(Safety)} for $k=1,\dots,m$, we have $A x_{k} \in S \quad \forall A \in U_{k-1}$,
    \item \textbf{(Learning)} the set of matrices $U_m$ is a singleton.
\end{enumerate}
\end{definition}


%

We now present our algorithm for checking the possibility of one-step safe learning (\algorithmref{alg:one-step}).
The proof of correctness of
\algorithmref{alg:one-step} is given in \theoremref{thm:one-step}.

\begin{remark}
As \theoremref{thm:one-step} will demonstrate, the particular choice of the parameter $\varepsilon\in(0,1]$ in the input to \algorithmref{alg:one-step} does not affect the detection of one-step safe learning by this algorithm. However, a smaller $\varepsilon$ leads to a lower cost of learning. Therefore, in practice, $\varepsilon$ should be chosen positive and as small as possible without causing the matrix $X$ in line~\ref{eq:x_mat} to be ill conditioned.
\end{remark}

%
%
%
\begin{algorithm2e}[htb]\label{alg:one-step}
\LinesNumbered
\SetKwInOut{Input}{Input}
\SetKwInOut{Output}{Output}
\SetKw{Break}{break}
\DontPrintSemicolon
\Input{polyhedra $S \subset \R^n$ and $U_A \subset \R^{n \times n}$, cost vector $c \in \R^n$, and a constant $\varepsilon \in (0, 1]$.}
\Output{A matrix $\trueA \in \R^{n \times n}$ or a declaration that one-step safe learning is impossible.}
\For{$k = 0,\dots,n-1$}{
$D_k \gets \{ (x_j, y_j) \mid j=1,\ldots,k\}$\;
$U_k \gets \{ A \in U_0 \mid  Ax_j = y_j, \quad j = 1,\dots,k\}$\;
\If{$U_k$ is a singleton (cf. \lemmaref{lem:uniqueness})}{
\Return the single element in $U_k$ as $\trueA$\; \label{alg:early_return}
}
Let $x_k^\star$ be the projection onto $x$-space of an optimal solution to problem \eqref{eq:one-step LP} with data $D_k$\;
\eIf{$x_k^\star$ is linearly independent from $\{x_1,\dots,x_k\}$}{
$x_{k+1} \gets x_k^\star$
}{
Let $S^1_k$ be the projection onto $x$-space of the feasible region of problem \eqref{eq:one-step LP} with data $D_k$\;
Compute a basis $B_k \subset S^1_k$ of $\spanof{S^1_k}$ (cf. \theoremref{thm:basis})\;
\For{$z_j \in B_k$}{
\If{$z_j$ is linearly independent from $\{x_1,\dots,x_k\}$}{
$x_{k+1} \gets (1-\varepsilon)  x_k^\star + \varepsilon z_j$\;
\Break
}
}
\If{no $z_j \in B_k$ is linearly independent from $\{x_1, \dots, x_k\}$}{
\Return one-step safe learning is impossible\; \label{alg:learning_impossible}
}
}
Observe $y_{k+1} \gets \trueA x_{k+1}$\;
}
Define matrix $X = [x_1,\dots,x_n]$ \label{eq:x_mat}\;
Define matrix $Y = [y_1,\dots,y_n]$\;
\Return $\trueA = YX^{-1}$ \label{alg:full_return}
\caption{One-Step Safe Learning Algorithm}
\end{algorithm2e}

\algorithmref{alg:one-step} invokes two subroutines which
we present next in \lemmaref{lem:uniqueness} and \theoremref{thm:basis}.

\begin{lemma}\label{lem:uniqueness}
Let $A \in \R^{m \times n}$, $B \in \R^{m \times p}$, $c \in \R^m$, and define the polyhedron
\[P \defn \{ x \in \R^n \mid  \exists y \in \R^p \quad \textnormal{s.t.} \quad A x + B y \leq c \}.\]
The problem of checking if $P$ is a singleton can be reduced to solving $2n$ linear programs.
\end{lemma}
\myproof{
For each $i=1,\dots,n$, maximize and minimize the $i$-th coordinate of $x$ over $P$.
It is straightforward to check that $P$ is a singleton if and only if the optimal values of these two linear programs coincide for every $i=1,\dots,n$.
}
%
%
For stating our next theorem, we use the following notation: given a set $P \subseteq \R^n$, let $\spanof{P}$ denote
the set of all linear combinations of points in $P$.
\begin{theorem}\label{thm:basis}
Let $A \in \R^{m \times n}$, $B \in \R^{m \times p}$, $c \in \R^m$, and define the polyhedron
\[P \defn \{ x \in \R^n \mid  \exists y \in \R^p \quad \textnormal{s.t.} \quad A x + B y \leq c \}.\]
One can find a basis of $\spanof{P}$ contained within $P$ by solving at most $2n^2$ linear programs.
\end{theorem}
\myproof{
We form the desired basis $\{e_i\}$ iteratively and with an inductive argument.
Let $e_1$ be any nonzero vector in $P$ (existence of such a vector can be checked by the argument in the proof of \lemmaref{lem:uniqueness}); if there is no such vector, we return the empty set.
Let $\{e_1,\dots,e_k\}$ be a linearly independent set in $P$.
We will either find an additional linearly independent vector $e_{k+1} \in P$, or show that the dimension of the span of $P$ is $k$.
Let $x$, $x^+$, and $x^-$ be variables in $\R^n$, $y^+$ and $y^-$ be variables in $\R^p$, and $\lambda^+$ and $\lambda^-$ be variables in $\R$.
Consider the following linear programming feasibility problem:
\begin{equation}\label{eq:one-step proof}
\begin{aligned}
& e_i^Tx = 0 \quad i=1,\dots,k\\
& x = x^+ - x^-\\
& Ax^+ + By^+ \leq \lambda^+ c \\
& Ax^- + By^- \leq \lambda^- c \\
& \lambda^+ \geq 0 \\
& \lambda^- \geq 0.
\end{aligned}
\end{equation}
Let $F\subseteq \R^n$ be the projection onto $x$-space of the feasible region of this problem.
We claim that $F = \{0\}$ if and only if the dimension of $\spanof{P}$ is $k$.
Moreover, if there is solution to \eqref{eq:one-step proof}
with $x \neq 0$,
then there is also a solution $(x, x^{\pm}, y^{\pm}, \lambda^{\pm})$ where $\lambda^+,\lambda^- \neq 0$.
In this case, either $\frac{x^+}{\lambda^+}$ or $\frac{x^+}{\lambda^+}$ can be taken as $e_{k+1}$.

Suppose first that the dimension $\spanof{P}$ is at least $k+1$; then there is a vector $\tilde{x}\in \spanof{P}$ that is linearly independent from $\{e_1,\dots,e_k\}$.
By subtracting the projection of $\tilde{x}$ onto $\spanof{\{e_1,\dots,e_k\}}$, we will find a nonzero vector $x \in \spanof{P}$ that is orthogonal to the vectors $e_1,\dots,e_k$.
We claim this vector $x$ is feasible to \eqref{eq:one-step proof} for some choice of $(x^\pm,y^\pm,\lambda^\pm)$.
Indeed, since $x \in \spanof{P}$, then
\[ x = \sum_{j=1}^r \lambda_j x_j, \]
for some vectors $x_1,\dots,x_r \in P$ and some nonzero scalars $\lambda_1,\dots,\lambda_r$.
For each $j$, as $x_j \in P$, there exists a vector $y_j$ such that $Ax_j + By_j \leq c$.
Let $J$ denote the set of indices $j$ such that $\lambda_j >0$.
It is easy to check that the assignment
\begin{equation}
\begin{aligned}
(x^+,y^+,\lambda^+) &= (\sum_{j \in J} \lambda_j x_j,\sum_{j \in J} \lambda_j y_j,\sum_{j \in J} \lambda_j) \:,\\
(x^-,y^-,\lambda^-) &= (- \sum_{j \notin J} \lambda_j x_j,- \sum_{j \notin J} \lambda_j y_j,- \sum_{j \notin J} \lambda_j)
\end{aligned} 
\end{equation}
satisfies system \eqref{eq:one-step proof}.
Hence, we have shown that if $F=\{0\}$ then the dimension of $\spanof{P}$ is $k$.
\par
To see the converse implication, suppose $x \neq 0$, and that the tuple $(x,x^\pm,y^\pm,\lambda^\pm)$ is feasible to system \eqref{eq:one-step proof}.
Without loss of generality we assume $\lambda^\pm\geq 1$; if not, we replace the tuple with
\begin{align}
    (x, x^\pm + \hat{x},y^\pm + \hat{y},\lambda^\pm + 1), \label{eq:make_lambda_nonzero}
\end{align}
where $\hat{x}$ and $\hat{y}$ are any vectors satisfying $A \hat{x} + B \hat{y} \leq c$.
Then, since $A \frac{x^+}{\lambda^+} + B \frac{y^+}{\lambda^+} \leq c$, the vector $\frac{x^+}{\lambda^+} \in P$.
By the same argument, $\frac{x^-}{\lambda^-} \in P$.
It follows from the orthogonality constraint of \eqref{eq:one-step proof} that at least one of the vectors $\frac{x^+}{\lambda^+}$ and $\frac{x^-}{\lambda^-}$ is linearly independent from $\{e_1,\dots,e_k\}$ and can be taken as $e_{k+1}$, also proving that the dimension of $\spanof{P}$ is at least $k+1$.
\par
Note that the condition $F=\{0\}$ can be checked by solving $2n$ linear programs (cf.\ the proof of \lemmaref{lem:uniqueness}); if $F\neq \{0\}$, then at least one of these $2n$ linear programs will return a tuple $(x,x^\pm,y^\pm,\lambda^\pm)$ where $x \neq 0$.
We then transform this tuple via \eqref{eq:make_lambda_nonzero}
to ensure that both $\lambda^+, \lambda^- \neq 0$
(we can take $\hat{x} = e_1$ and $\hat{y}$ to be any vector
such that $A e_1 + B \hat{y} \leq c$).
Since we cannot have more than $n$ linearly independent vectors in $\spanof{P}$, this procedure is repeated at most $n$ times.
}


Our next theorem is the main result of the section.
\begin{theorem}
\label{thm:one-step}
Given a safety region $S \subset \R^n$ and an uncertainty set $U_0 \subset \R^{n \times n}$, one-step safe learning is possible if and only if \algorithmref{alg:one-step} (with an arbitrary choice of $c \in \R^n$ and $\varepsilon \in (0, 1]$)
returns a matrix.
\end{theorem}
\myproof{
[``If'']
By construction, the sequence of measurements chosen by \algorithmref{alg:one-step} satisfies the first condition of \definitionref{def:one-step}, since the vectors $x_k^\star$ and $z_j$ are both contained in $S^1_k$ and any vector in $S^1_k$ will remain in the safety region under the action of all matrices in $U_k$; i.e. all matrices in $U_A$ that are consistent with the measurements made so far.
If \algorithmref{alg:one-step} terminates early at line~\ref{alg:early_return} for some iteration $k$, then clearly the uncertainty set $U_k$ is a singleton.
On the other hand, 
if we reach line~\ref{alg:full_return}, then we must have made $n$
linearly independent measurements $\{x_1,\dots,x_{n}\}$.
From this, it is clear that the set $\{ A \in U_0 \mid A x_j = y_j, j = 1, \dots, n \} = \{A_\star\}$.
\par
[``Only if''] Suppose \algorithmref{alg:one-step} chooses points $\{x_1,\dots,x_{m}\}$ where $m < n$ and terminates at line~\ref{alg:learning_impossible}.
Then it is clear from the algorithm that $\{x_1,\dots,x_{m}\}$ must form a basis of $\spanof{S^1_m}$ and that $U_m$ is not a singleton.
Take $\tilde{m}$ to be any nonnegative integer and $\{\tilde{x}_1,\dots,\tilde{x}_{\tilde{m}}\}$ to be any sequence that satisfies the first condition of \definitionref{def:one-step}.
For $k=1, \ldots, \tilde{m}$, let
\begin{align*}
    \tilde{U}_k &= \{ A \in U_0 \mid  A\tilde{x}_j = \trueA \tilde{x}_{j},  j = 1,\dots,k\} \:, \\
    \tilde{S}^1_k &= \{x \in S \mid Ax \in S,  \forall A\in \tilde{U}_k\} \:.
\end{align*}
First we claim that $\tilde{x}_k \in S^1_m$ for $k=1,\dots,\tilde{m}$.
We show this by induction.
It is clear that $\tilde{x}_1 \in S^1_m$ since $\tilde{x}_1 \in S^1_0$ and $S^1_0 \subseteq S^1_m$.
Now we assume $\tilde{x}_1,\dots,\tilde{x}_{k} \in S^1_m$ and show that $\tilde{x}_{k+1} \in S^1_m$.
Since $\{x_1,\dots,x_{m}\}$ forms a basis of $\spanof{S^1_m}$, it follows that for any matrix $A$, $Ax_j=\trueA x_j$ for $j=1,\dots,m$ implies $Ax = \trueA x$ for all $x \in S^1_m$.
In particular, for any matrix $A$, $Ax_j=\trueA x_j$ for $j=1,\dots,m$ implies $A\tilde{x}_{j} = \trueA \tilde{x}_{j}$ for all $j=1,\dots,k$.
It follows that $U_m \subseteq \tilde{U}_k$ and therefore, $\tilde{S}^1_k \subseteq S^1_m$.
By the first condition of \definitionref{def:one-step}, we must have $\tilde{x}_{k+1} \in \tilde{S}^1_k$, and thus, $\tilde{x}_{k+1} \in S^1_m$.
This completes the inductive argument and shows that $\tilde{x}_k \in S^1_m$ for $k=1,\dots,\tilde{m}$.
From this, it follows that $U_m \subseteq \tilde{U}_{\tilde{m}}$.
Recall that $U_m$ is not a singleton, thus $\tilde{U}_{\tilde{m}}$ is not a singleton either.
Therefore, the sequence $\{\tilde{x}_1,\dots,\tilde{x}_{\tilde{m}}\}$ does not satisfy the second condition of \definitionref{def:one-step}.
}

\begin{corollary}
Given a safety region $S \subset \R^n$ and an uncertainty set $U_0 \subset \R^{n \times n}$, if one-step safe learning is possible, then it is possible with at most $n$ measurements. 
\end{corollary}
\subsection{The Value of Exploiting Information on the Fly}
\label{sec:cost of learning}

In addition to detecting the possibility of safe learning, \algorithmref{alg:one-step}
attempts to minimize the overall cost of learning (i.e.,\ $\sum_{k=1}^{m} c^T x_k$) by exploiting 
information gathered at every step.
In order to demonstrate the value of using information online,
we construct an offline algorithm which 
chooses $n$ measurement vectors $x_1, \dots, x_n$ ahead of time based solely on $U_0$ and $S$,
and succeeds under the assumption that $S_0^1$ contains a basis of $\R^n$.
%
%
%
%
\begin{algorithm2e}[ht]\label{alg:offline-one-step}
\LinesNumbered
\SetKwInOut{Input}{Input}
\SetKwInOut{Output}{Output}
\SetKw{Break}{break}
\DontPrintSemicolon
\Input{polyhedra $S \subset \R^n$ and $U_0 \subset \R^{n \times n}$, cost vector $c \in \R^n$, and a constant $\varepsilon \in (0, 1]$.}
\Output{A matrix $\trueA \in \R^{n \times n}$ or failure.}
\If{$S_0^1$ does not contain a basis of $\R^n$ (cf. \theoremref{thm:basis})}{
\Return failure
}
Compute a basis $\{ z_1, \dots, z_n \} \subset S_0^1$ of $\R^n$\;
Let $x_0^\star$ be the projection onto $x$-space of an optimal solution to problem \eqref{eq:one-step LP} with data $D_0$\;
Set $x_k = (1-\varepsilon) x_0^\star + \varepsilon z_k$ for $k=1, \dots, n$\;
Observe $y_k \gets A_\star x_k$ for $k=1, \dots, n$\;
Define matrix $X = [x_1,\dots,x_n]$\;
Define matrix $Y = [y_1,\dots,y_n]$\;
\Return $\trueA = YX^{-1}$
\caption{Offline Safe Learning Algorithm}
\end{algorithm2e}

%
As $\varepsilon$ tends to zero, the cost of Algorithm~\ref{alg:offline-one-step}
approaches $n c^T x_0^\star$, where $x_0^\star$ is a minimum cost measurement vector
in $S_0^1$.
Therefore, $n c^T x_0^\star$
serves as an \emph{upper bound} on the cost incurred by \algorithmref{alg:one-step}.
We note that $n c^T x_0^\star$ is also the minimum cost achievable by any one-step safe offline algorithm
that takes $n$ measurements, 
since all measurement vectors $\{x_k\}$ of such an algorithm must come from $S^1_0$.

By assuming that we know $\trueA$,
we can also compute a \emph{lower bound} on the cost of one-step safe learning of any algorithm
that takes $n$ measurements. 
Let $S^1(\trueA) = \{ x \in S \mid \trueA x \in S \}$ be the true one-step safety region of $\trueA$.
Let $x^\star$ be an optimal solution to the linear program that minimizes $c^T x$ over $S^1(\trueA)$.
Then, clearly, if we must pick $n$ points that are all one-step safe, we cannot 
achieve cost lower than $n c^T x^\star$.

\subsection{Numerical Example}\label{sec:one-step example}
\begin{figure}
\figureconts
{fig:one-step}
{\caption{One-step safe learning associated with the numerical example in \sectionref{sec:one-step example}.}}
{%
\subfigure[$S^1_k$ grows with $k$.][c]{%
\label{fig:one-step S}
\includegraphics[width=.5\textwidth -.5em]{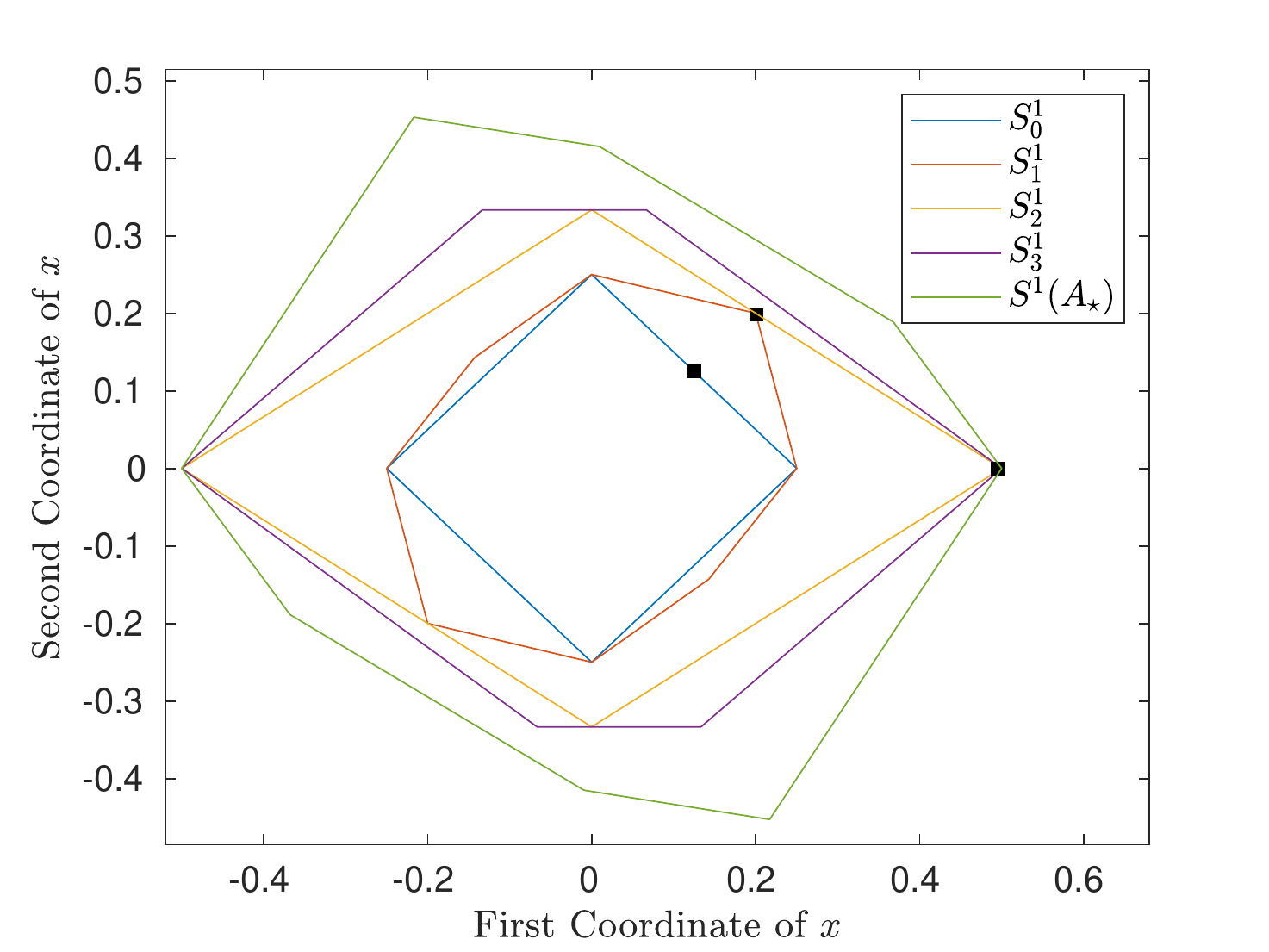}
} 
\subfigure[$U_k$ shrinks with $k$.][c]{%
\label{fig:one-step U}
\includegraphics[width=.5\textwidth -.5em]{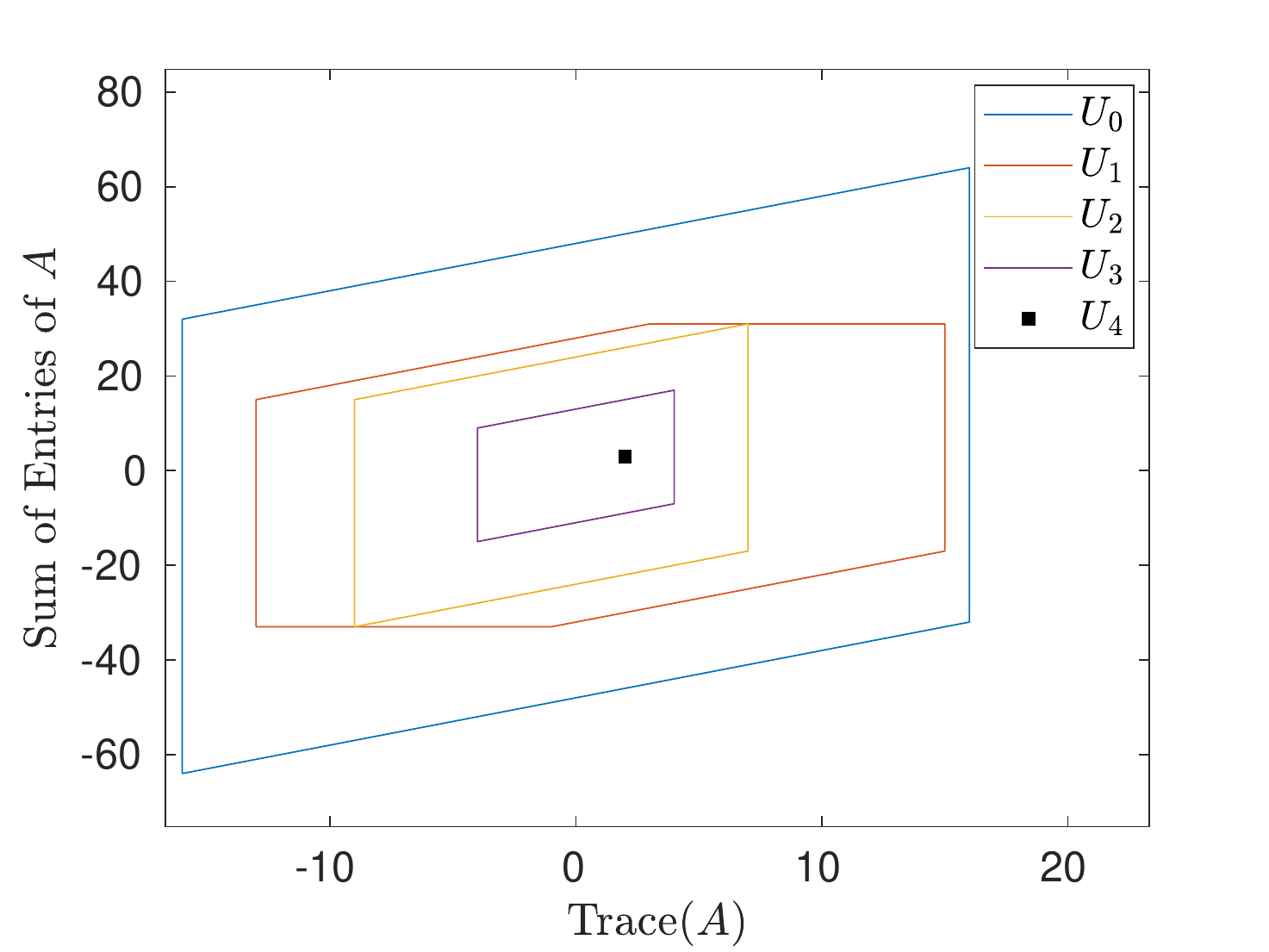}
}
}
\end{figure}

We present a numerical example with $n=4$.
Here, we take $U_0 = \{ A \in \R^{4 \times 4} \mid |A_{ij}| \leq 4 \quad \forall i,j \}$,
$S = \{ x \in \R^4 \mid \norm{x}_{\infty} \leq 1 \}$,
and $c=(-1,-1,0,0)^T$.
We choose the matrix $\trueA$ uniformly at random among integer matrices in $U_0$:
\[
\trueA = \begin{bmatrix*}[r]
 2 &  1 &  4 &  2\\
 2 & -3 & -1 & -2\\
-2 & -3 &  1 &  0\\
 2 &  0 & -2 &  2
\end{bmatrix*}.
\]
In this example, \algorithmref{alg:one-step}
takes four steps to safely recover $\trueA$.
The projection onto the first two dimensions of the four vectors that \algorithmref{alg:one-step}
selects are plotted in \figureref{fig:one-step S} (note that two of the points are very close to each other).
Because of the cost vector $c$, points higher and further to the right in the plot have lower measurement cost.
Also plotted in \figureref{fig:one-step S} are the projections onto the first two dimensions
of the sets $S^1_k$ for $k\in\{0, 1, 2, 3\}$
and of the set $S^1(\trueA)$,
the true one-step safety region of $\trueA$.
In \figureref{fig:one-step U}, we plot $U_k$ (the remaining uncertainty after making $k$ measurements) for $k\in\{0, 1, 2, 3, 4\}$; we draw a two-dimensional projection of these sets of matrices by looking at the trace and the sum of the entries of each matrix in the set.
Note that $U_4$ is a single point since we have recovered the true dynamics after the fourth measurement.

The cost of learning (i.e.,\ $\sum_{i=1}^{4} c_i^T x_i$) 
is $-1.0000$ for the offline algorithm (\algorithmref{alg:offline-one-step}), and $-1.6385$ for \algorithmref{alg:one-step}.
The lower bound on the cost of learning is $-2.2264$ (cf.\ \sectionref{sec:cost of learning}).
We can see that the value of exploiting information on the fly is significant.
\section{Two-Step Safe Learning of Linear Systems}
In this section, we again focus on learning the linear dynamics in \eqref{eq:linear dynamics}.
However, unlike the previous section, we are interested in making queries to the system that are two-step safe.
The advantage of this formulation is that we may have fewer system resets and can potentially learn the dynamics with lower measurement cost.

More formally, in the two-step safe learning problem, we have as input a polyhedral safety region $S \subset \R^n$ given in the form of \eqref{eq:S polyhedron}, an objective function representing measurement cost which for simplicity we again take to be a linear function $c^T x$, and an uncertainty set $U_0 \subset \R^{n \times n}$ to which we assume $\trueA$ belongs.
In this section, we take $U_0$ to be an ellipsoid; this means that there is a strictly convex quadratic function $q: \R^{n \times n} \rightarrow \R$ such that:
\[ U_0 = \left \{ A \in \R^{n \times n} \mid  q(A) \leq 0 \right \}.\]
An example of such an uncertainty set is $U_0 =  \{ A \in \R^{n \times n} \mid \norm{A-A_0}_F \leq \gamma\} $, where $A_0$ is a nominal matrix, $\gamma$ is a positive scalar, and $\norm{\cdot}_F$ refers to the Frobenius norm.
Having collected $k$ safe length-two trajectories $\{ (x_j, \trueA x_j, \trueA^2 x_j) \}_{j=1}^{k}$, 
our uncertainty around $\trueA$ reduces to:
\begin{align}
    U_k = \{ A \in U_0 \mid A x_j = \trueA x_j, A^2 x_j = \trueA^2 x_j, j=1, \dots, k \}. \label{eq:U_k_two_step}
\end{align}
The optimization problem we would like to solve to find the next best two-step safe query point is the following:
\begin{equation}\label{two-step}
\begin{aligned}
\min_{x} \quad & c^T x\\
\textrm{s.t.} \quad & x \in S \\
& A x \in S \quad \forall A \in U_k \\
& A^2 x \in S \quad \forall A \in U_k.
\end{aligned}
\end{equation}

\subsection{Reformulation via the S-Lemma}
In this subsection, we derive a tractable reformulation of problem \eqref{two-step}.
\begin{theorem}\label{thm:two-step}
Problem \eqref{two-step} can be reformulated as a semidefinite program.
\end{theorem}
Our proof makes use to the S-lemma 
\citep[see e.g.,][]{s_lemma} which we recall next.
\begin{lemma}[S-lemma]
For two quadratics functions $q_a$ and $q_b$, if there exists a point $\bar{x}$ such that $q_a(\bar{x}) < 0$, then the implication
\[ \forall x, \left[ q_a(x) \leq 0 \Rightarrow q_b(x) \leq 0 \right] \]
holds if and only if there exists a scalar $\lambda \geq 0$ such that
\[\lambda q_a(x) - q_b(x) \geq 0 \quad \forall x .\]
\end{lemma}
\myproofof{of \theoremref{thm:two-step}}{
Note that the set of equations
\[Ax_j = \trueA x_j, \quad  A^2 x_j = \trueA^2 x_j \quad j =1,\dots,k\]
in the definition of $U_k$ in \eqref{eq:U_k_two_step} is equivalent to the set of linear equations 
\begin{equation}\label{eq:two-step linear constraints}
    Ax_j = \trueA x_j, \quad  A (\trueA x_j) = \trueA^2 x_j \quad j =1,\dots,k.
\end{equation}
If there is only one matrix in $U_0$ that satisfies all of the equality constraints in \eqref{eq:two-step linear constraints} (a condition that can be checked via a simple modification of \lemmaref{lem:uniqueness}), then we have found $\trueA$ and \eqref{two-step} becomes a linear program.
Therefore, let us assume that more than one matrix in $U_0$ satisfies the constraints in \eqref{eq:two-step linear constraints}.
In order to apply the S-lemma, we need to remove these equality constraints, a task that we accomplish via variable elimination.
Let $\hat{n}$ be the dimension of the affine subspace of matrices that satisfy the constraints in \eqref{eq:two-step linear constraints} and let $\hat{A} \in \R^{n \times n}$ be an arbitrary member of this affine subspace.
Let $A_1,\dots,A_{\hat{n}}  \in \R^{n \times n}$ be a basis of the subspace
\[\{ A \in \R^{n \times n} \mid Ax_j = 0,\quad A(\trueA x_j) = 0 \quad j=1,\dots,k \}. \]
Consider an affine function $g : \R^{\hat{n}} \rightarrow \R^{n \times n}$ defined as follows:
\[
g(\hat{a}) \defn \hat{A} + \sum_{i=1}^{\hat{n}} \hat{a}_i A_i.
\]
The function $g$ has the properties that it is injective and that for each $A$ that satisfies the equality constraints, there must be a vector $\hat{a}$ such that $A = g(\hat{a})$.
In other words, the function $g$ is simply parametrizing the affine subspace of matrices that satisfy the equality constraints.
Now we can reformulate \eqref{two-step} as:
\begin{equation}\label{eq:two-step parametrization}
\begin{aligned}
\min_{x} \quad & c^T x\\
\textrm{s.t.} \quad & x \in S \\
& g(\hat{a}) x \in S \quad \forall \hat{a} \quad \textrm{s.t.} \quad q ( g(\hat{a})) \leq 0 \\
& g(\hat{a})^2 x \in S \quad \forall \hat{a} \quad \textrm{s.t.} \quad q ( g(\hat{a})) \leq 0.
\end{aligned}
\end{equation}
Let $\hat{q} \defn q \circ g$.
Since $q$ is a strictly convex quadratic function and $g$ is an injective affine map, $\hat{q}$ is also a strictly convex quadratic function.
Since we are under the assumption that there are multiple matrices in $U_0$ that satisfy the equality constraints, there must be a vector $\bar{a} \in \R^{\hat{n}}$ such that $\hat{q}(\bar{a}) < 0$.
To see this, take $\bar{a}_1 \neq \bar{a}_2$ such that $\hat{q}(\bar{a}_1),\hat{q}(\bar{a}_2) \leq 0$.
It follows from strict convexity of $\hat{q}$ that $\hat{q}(\frac{1}{2} (\bar{a}_1 + \bar{a}_2)) < 0$.
Using the definition of $S$, problem \eqref{eq:two-step parametrization} can be rewritten as:
\begin{equation}\label{eq:two-step bilevel}
\begin{aligned}
\min_{x} \quad & c^T x\\
\textrm{s.t.} \quad & h_i^T x \leq b_i \quad i = 1, \dots ,r \\
& \begin{bmatrix} \max_{\hat{a}} \quad & h_i^T g(\hat{a}) x  \\ 
\textrm{s.t.} \quad & \hat{q}(\hat{a}) \leq 0 \end{bmatrix} \leq b_i \quad i = 1, \dots ,r \\
& \begin{bmatrix} \max_{\hat{a}} \quad & h_i^T g(\hat{a})^2 x  \\ 
\textrm{s.t.} \quad & \hat{q}(\hat{a}) \leq 0 \end{bmatrix} \leq b_i \quad i = 1, \dots ,r.
\end{aligned}
\end{equation}
Let $q_{1,i}(\hat{a};x) = h_i^T g(\hat{a}) x - b_i$ and $q_{2,i}(\hat{a};x) = h_i^T g(\hat{a})^2 x - b_i$.
We consider these functions as quadratic functions of $\hat{a}$ parametrized by $x$.
Note that the coefficients of $q_{1,i}$ and $q_{2,i}$ depend affinely on $x$.
Using logical implications, problem \eqref{eq:two-step bilevel} can be rewritten as:
\begin{equation}\label{eq:two-step implications}
\begin{aligned}
\min_{x} \quad & c^T x\\
\textrm{s.t.} \quad & h_i^T x \leq b_i \quad i = 1, \dots ,r \\
& \forall \hat{a}, \left[ \hat{q}(\hat{a}) \leq 0 \Rightarrow q_{1,i}(\hat{a};x) \leq 0 \right] \quad i = 1, \dots ,r \\
& \forall \hat{a}, \left[ \hat{q}(\hat{a}) \leq 0 \Rightarrow q_{2,i}(\hat{a};x) \leq 0 \right] \quad i = 1, \dots ,r.
\end{aligned}
\end{equation}
Now we use the S-lemma to reformulate an implication between quadratic inequalities as a constraint on the global nonnegativity of a quadratic function.
Note that as we have already argued for the existence of a vector $\bar{a}$ such that $\hat{q}(\bar{a}) < 0$, the condition of the S-lemma is satisfied.
After introducing variables $\lambda_{1,i}$ and $\lambda_{2,i}$  for $i=1,\dots,r$, we apply the S-lemma $2r$ times to reformulate \eqref{eq:two-step implications} as the following program:
\begin{equation}\label{eq:two-step sdp}
\begin{aligned}
\min_{x,\lambda} \quad & c^T x\\
\textrm{s.t.} \quad & h_i^T x \leq b_i \quad i = 1, \dots ,r \\
& \lambda_{1,i} \hat{q}(\hat{a}) - q_{1,i}(\hat{a};x) \geq 0 \quad \forall \hat{a} \quad i = 1, \dots ,r \\
& \lambda_{2,i} \hat{q}(\hat{a}) - q_{2,i}(\hat{a};x) \geq 0 \quad \forall \hat{a} \quad i = 1, \dots ,r \\
& \lambda_{1,i} \geq 0, \quad \lambda_{2,i} \geq 0 \quad i=1,\dots,r.
\end{aligned}
\end{equation}
It is a standard procedure to convert the constraint that a quadratic function is globally nonnegative into a semidefinite constraint.
Note that the coefficients of $q_{1,i}$ and $q_{2,i}$ depend affinely on $x$; this results in linear matrix inequalities when \eqref{eq:two-step sdp} is converted into a semidefinite program.
}

\subsection{Numerical Example}
\label{sec:two-step example}

\begin{figure}[t]
\figureconts
{fig:two-step}
{\caption{Two-step safe learning associated with the numerical example in \sectionref{sec:two-step example}.}}
{%
\subfigure[$S^2_k$ grows with $k$.][c]{%
\label{fig:two-step S}
\includegraphics[width=.5\textwidth -.5em]{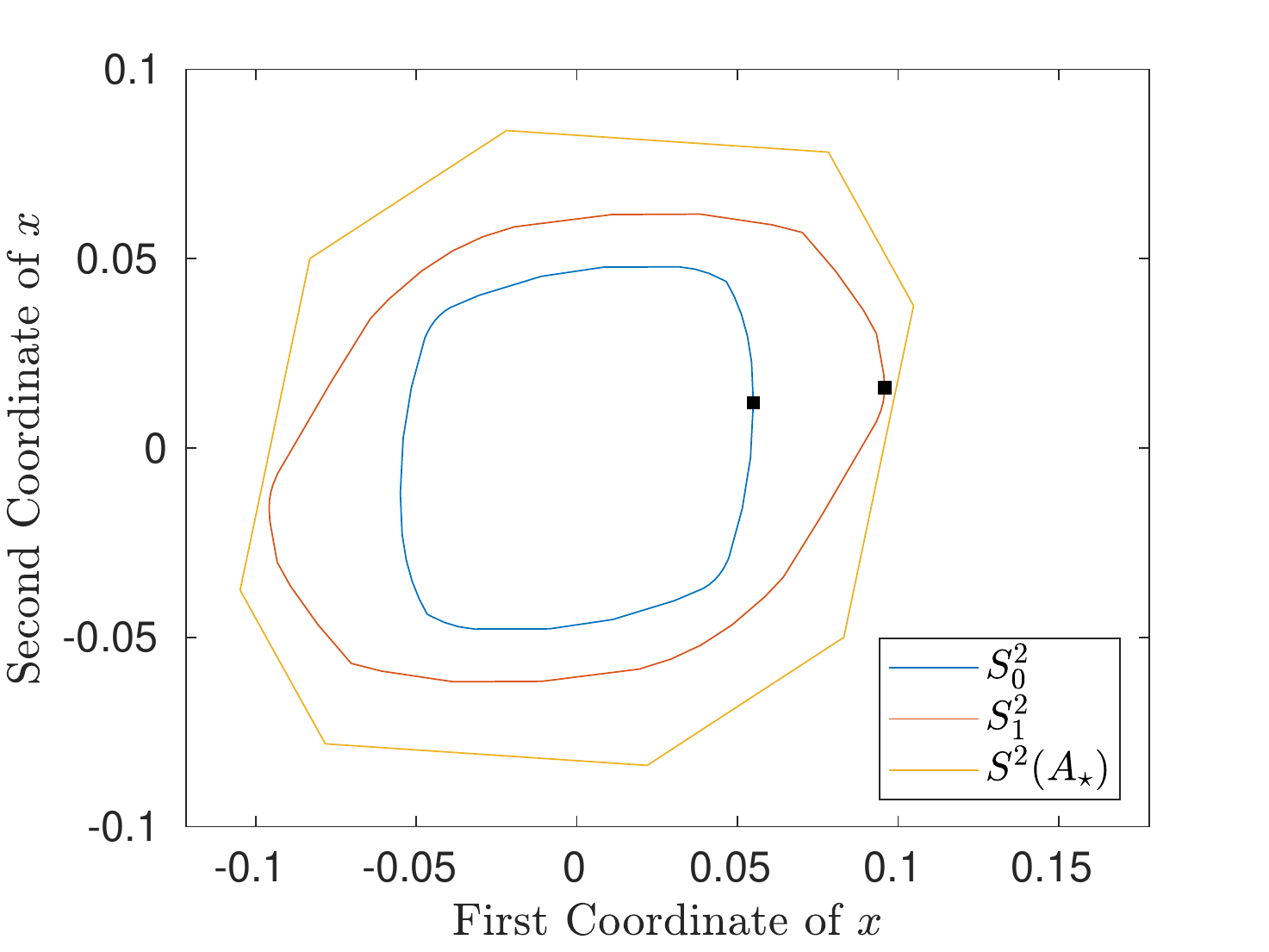}
} 
\subfigure[$U_k$ shrinks with $k$.][c]{%
\label{fig:two-step U}
\includegraphics[width=.5\textwidth -.5em]{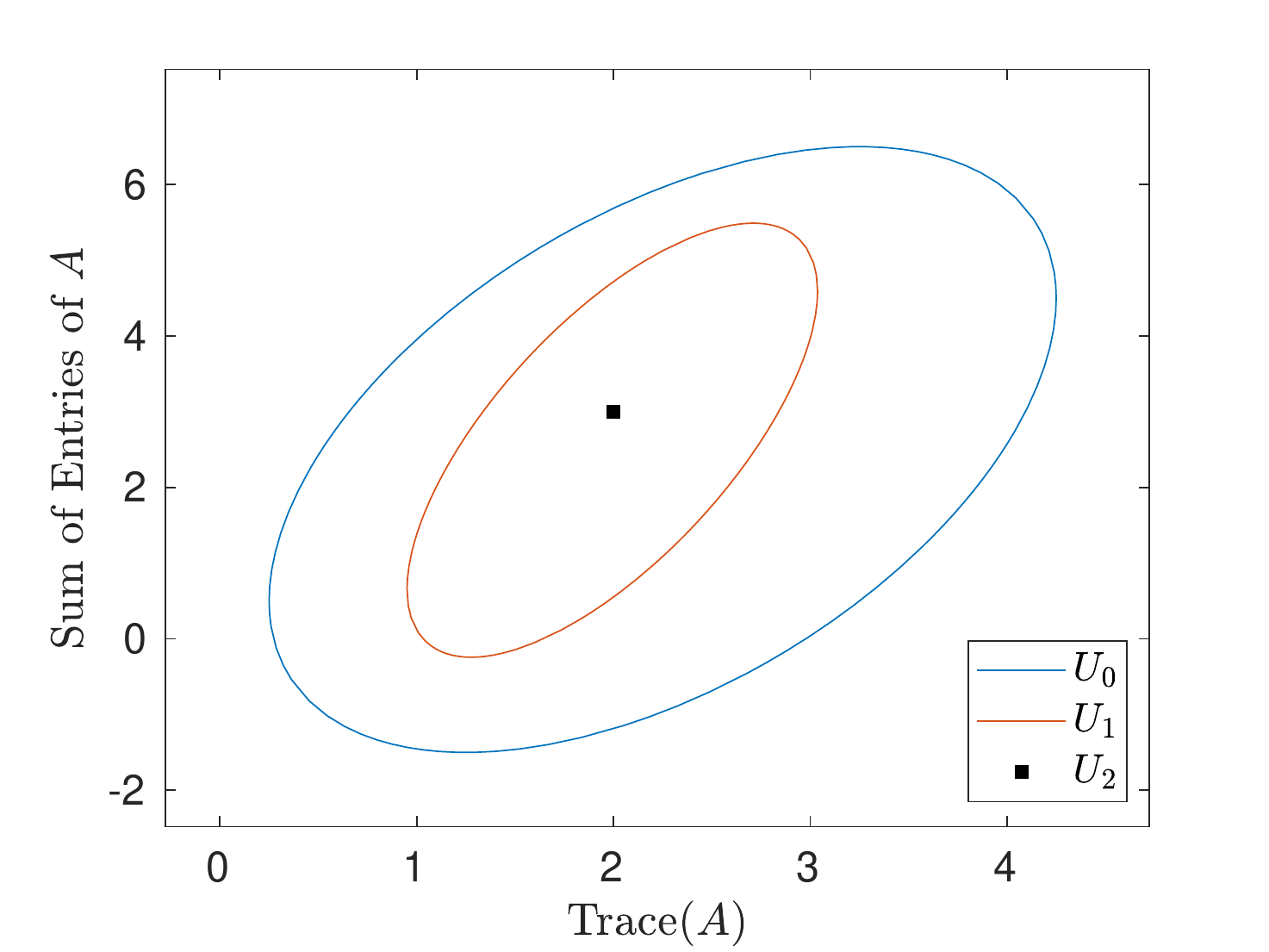}
}
}
\end{figure}

We present a numerical example, again with $n=4$.
Here we choose a nominal matrix
\[
A_0 = \begin{bmatrix*}[r]
 2.25 &  0.75 &  4.25 &  1.75\\
 2.25 & -3.25 & -1.25 & -2.25\\
-2.00 & -2.75 &  1.25 &  0.00\\
 1.75 & -0.25 & -2.00 &  2.00
\end{bmatrix*}
\]
and let $U_0 = \{ A \in \R^{4 \times 4} \mid \norm{A-A_0}_F \leq 1 \}$.
We let $S = \{ x \in \R^4 \mid |x_i| \leq 1, i=1,\dots,4 \}$ and $c=(-1,0,0,0)^T$.
We choose the true matrix $\trueA$ to be the same matrix used in \sectionref{sec:one-step example}.
\par
In this example, by solving two semidefinite programs, we learn the true matrix $\trueA$ by making two measurements that are each two-step safe.
In other words, we choose $x_1\in\R^4$, observe $\trueA x_1$, $\trueA^2 x_1,$ and then choose $x_2\in\R^4$, and observe $\trueA x_2$ and $\trueA^2 x_2$.
We can verify that we have recovered $\trueA$ if $\{ x_1, \trueA x_1 , x_2, \trueA x_2\}$ are all linearly independent, which is the case.
The projection onto the first two dimensions of the two measurements $x_1$ and $x_2$ that our semidefinite programs choose are plotted in \figureref{fig:two-step S}.
Because of the cost vector $c$, points further to the right in the plot have lower measurement cost.
Also plotted are the projections onto the first two dimensions of the sets:
\[ \begin{aligned}
S^2_0 &= \{ x \in S \mid A x \in S, A^2 x \in S \quad \forall A \in U_0\}, \\
S^2_1 &= \{ x \in S \mid A x \in S, A^2 x \in S \quad \forall A \in U_1 \}, \\
S^2(\trueA) &= \{ x \in S \mid \trueA x \in S, \trueA^2 x \in S \}.
\end{aligned} \]
The first two sets are the projections onto $x$-space of the feasible regions of our two semidefinite programs.
The third set is the true two-step safety region of $\trueA$.
In \figureref{fig:two-step U}, we plot $U_k$ (the remaining uncertainty after observing $k$ trajectories of length two) for $k\in\{0, 1, 2\}$; we draw a two-dimensional projection of these sets of matrices by looking at the trace and the sum of the entries of each matrix in the set.
Note that $U_2$ is a single point since we have recovered the true dynamics after observing the second trajectory.
The cost of learning (i.e.,\ $c^T x_1 + c^T x_2$) 
is $-0.1508$.

We can construct an analogue of the offline \algorithmref{alg:offline-one-step} by only making measurements from $S^2_0$.
This approach would first pick the optimal point in $S^2_0$ (i.e., $x_1$), and then another vector in $S^2_0$ close to $x_1$, but linearly independent from it.
The cost of learning for this offline approach would be $2c^T x_1 = -0.1099$.
Finally, we can again find a lower bound on the cost of learning
of any algorithm that makes two measurements (that are each two-step safe)
by assuming we know $\trueA$ ahead of time and optimizing $c^T x$ over $S^2(\trueA)$; in this example, the lower bound is $-0.2097$.
Here, again, we see that by using information on the fly, we can succeed at safe learning at a lower cost than the offline approach.
\section{One-Step Safe Learning of Nonlinear Systems}
\label{sec:nonlinear}

We consider the problem of safely learning a dynamical system
of the form $x_{t+1} = \truef(x_t)$, where
\begin{align}
  \truef(x) = \trueA x + \trueg(x), \label{eq:nonlinear_model}
\end{align}
for some matrix $\trueA \in \R^{n \times n}$
and some possibly nonlinear map $\trueg :  \R^n \rightarrow \R^n$.
We take our safety region $S \subset \R^n$ to be the same as \eqref{eq:S polyhedron}.
Our initial knowledge about $\trueA, \trueg$ is membership in the sets
\begin{align*}
    U_{0,A} &\defn \left \{ A \in \R^{n \times n} \mid  \Tr(V_j^T A) \leq v_j \quad j = 1, \dots ,s \right \}, \\
    U_{0,g} &\defn \{ g : \R^n \rightarrow \R^n \mid \norm{g(x)}_{\infty} \leq \gamma \norm{x}_p^d \quad \forall x \in S \}.
\end{align*}
Here, $p \geq 1$ is either $+\infty$ or a rational number,
$\gamma$ is a given positive constant, and $d$ is a given nonnegative integer.
The use of the $\norm{\cdot}_\infty$ on $g$ in the definition of $U_{0, g}$
simplifies some of the following analysis,
though an extension to other semidefinite representable norms is possible.
Note that by taking $d=0$ e.g.,\ our model of uncertainty captures any
map $f$ which is bounded on $S$.

%
%


Again for simplicity, we assume a linear measurement cost $c^T x$ for some vector $c \in \R^n$.
Having collected $k$ safe measurements $\{(x_j, y_j)\}_{j=1}^{k}$
with $y_j = \truef(x_j)$,
the optimization problem we are interested in solving to find the next cheapest one-step safe measurement is:
\begin{equation}\label{eq:nonlinear}
\begin{aligned}
\min_{x} \quad & c^T x &\\
\textrm{s.t.} \quad & x \in S &\\
& A x + g(x) \in S \quad \forall~(A,g) \in \{ A \in U_{0, A}, g \in U_{0, g} \mid Ax_j + g(x_j) = y_j \quad j = 1,\dots,k \}.
\end{aligned}
\end{equation}
\subsection{Reformulation as a Second-Order Cone Program}
Our main result of this section is to derive a tractable reformulation of problem~\eqref{eq:nonlinear}.
\begin{theorem}\label{thm:nonlinear socp}
Problem \eqref{eq:nonlinear} can be reformulated as a second-order cone program.
\end{theorem}
\myproof{
We start by rewriting problem~\eqref{eq:nonlinear} using the definition of $S$:
\begin{equation}\label{eq:nonlinear bilevel}
\begin{aligned}
\min_{x} \quad & c^T x\\
\textrm{s.t.} \quad & h_i^T x \leq b_i \quad i = 1, \dots ,r \\
& \begin{bmatrix} \max_{A,g} \quad & h_i^T (A x + g(x))  \\ 
\textrm{s.t.} \quad & \Tr(V_j^T A) \leq v_j \quad j = 1, \dots ,s \\
& \norm{g(x)}_{\infty} \leq \gamma \norm{x}_p^d \quad \forall x \in S \\
& Ax_k + g(x_k) = y_k \quad k = 1,\dots,m 
\end{bmatrix} \leq b_i \quad i = 1, \dots ,r.
\end{aligned}
\end{equation}

Note that in the inner maximization problem in \eqref{eq:nonlinear bilevel},
the variable $x$ is fixed.
We claim that if $x \not\in \{x_1, \dots, x_k\}$, then
\begin{align}
&\begin{bmatrix} \max_{A,g} \quad & h_i^T (A x + g(x))  \\ 
\textrm{s.t.} \quad & \Tr(V_j^T A) \leq v_j \quad j = 1, \dots ,s \\
& \norm{g(x)}_{\infty} \leq \gamma \norm{x}_p^d \quad \forall x \in S \\
& Ax_k + g(x_k) = y_k \quad k = 1,\dots,m 
\end{bmatrix} \label{eq:three_brackets} \\ 
&\qquad = \begin{bmatrix} \max_{A,g} & h_i^T A x  \\ 
\textrm{s.t.} & \Tr(V_j^T A) \leq v_j \quad \forall j \\
& Ax_k + g(x_k) = y_k \quad \forall k \\
& \norm{g(x)}_{\infty} \leq \gamma \norm{x}_p^d \quad \forall x \in S
\end{bmatrix} + \begin{bmatrix}
\max_{A,g} & h_i^T g(x)  \\ 
\textrm{s.t.} & \Tr(V_j^T A) \leq v_j \quad \forall j \\
& Ax_k + g(x_k) = y_k \quad \forall k \\
& \norm{g(x)}_{\infty} \leq \gamma \norm{x}_p^d \quad \forall x \in S
\end{bmatrix}. \nonumber 
\end{align}
It is clear that the left-hand side is upper bounded by the right-hand side.
To show the reverse inequality, let
$(A_1, g_1)$ (resp.\ $(A_2, g_2)$) be feasible to the first (resp. second) problem on the right-hand side (if any of these of these problems is infeasible, then the inequality we are after is immediate).
Now let 
\begin{align*}
    \hat{g}_2(x) = \begin{cases}
        g_2(x) &\text{if } x \not\in \{x_1, \dots, x_k\} \\
        y_k - A_1 x_k &\text{if } x = x_k.
    \end{cases}
\end{align*}
It is straightforward to check that the pair $(A_1, \hat{g}_2)$ is feasible to the left-hand side of \eqref{eq:three_brackets},
therefore proving \eqref{eq:three_brackets}.

We now focus on reformulating each term on the right-hand side of \eqref{eq:three_brackets}, again under the assumption that $x \not\in \{ x_1,\dots,x_m\}$.
Using the constraint on $g$, the first term can be rewritten as follows:
\begin{equation}\label{eq:max A nonlinear}
\begin{aligned}
\max_{A} \quad & h_i^T A x  \\ 
\textrm{s.t.} \quad & \Tr(V_j^T A) \leq v_j \quad j = 1, \dots ,s \\
& \norm{Ax_k - y_k}_{\infty} \leq \gamma \norm{x_k}_p^d \quad k = 1,\dots,m .
\end{aligned}
\end{equation}
Note that \eqref{eq:max A nonlinear} is a linear program as it is equivalent to:
\begin{equation}\label{eq:nonlinear primal}
\begin{aligned}
\max_{A} \quad & h_i^T A x  \\ 
\textrm{s.t.} \quad & \Tr(V_j^T A) \leq v_j \quad j = 1, \dots ,s \\
& (Ax_k - y_k)_{l} \leq \gamma  \norm{x_k}_p^d \quad k = 1,\dots,m \quad l = 1,\dots,n\\
& -(Ax_k - y_k)_{l} \leq \gamma  \norm{x_k}_p^d \quad k = 1,\dots,m \quad l = 1,\dots,n.
\end{aligned}
\end{equation}
Here, the notation $(Ax_k - y_k)_{l}$ represents the $l$-th coordinate of the vector $(Ax_k - y_k)$.
Following the same approach as in \sectionref{sec:one-step}, we proceed by taking the dual of this linear program.
For $j=1,\dots,s$, $k=1,\dots,m$, and $l=1,\dots,n$, let $\mu_j,\eta^+_{kl},\eta^-_{kl} \in \R$ be dual variables.
The dual of problem \eqref{eq:nonlinear primal} reads:
\begin{equation}
\begin{aligned}
\min_{\mu,\eta^+,\eta^-} \quad & \sum_j \mu_j v_j + \sum_{kl} \eta^+_{kl} (\gamma\norm{x_k}_p^d + (y_k)_l )  + \sum_{kl} \eta^-_{kl} (\gamma\norm{x_k}_p^d - (y_k)_l )\\ 
\textrm{s.t.} \quad & x h_i^T = \sum_j \mu_j V_j^T + \sum_{kl} \eta^+_{kl} x_k e_l^T - \sum_{kl} \eta^+_{kl} x_k e_l^T \\
& \mu \geq 0, \quad \eta^+ \geq 0, \quad \eta^- \geq 0,
\end{aligned}
\end{equation}
where $e_l$ is the $l$-th coordinate vector.
Now we turn our attention to the second term on the right-hand side of \eqref{eq:three_brackets}.
After eliminating the irrelevant constraints, the problem can be rewritten as:
\begin{equation}
\begin{aligned}
\max_{g} \quad &  h_i^T g(x)  \\ 
\textrm{s.t.} \quad & \norm{g(x)}_{\infty} \leq \gamma \norm{x}_p^d.
\end{aligned}
\end{equation}
Recall that the dual norm of $\norm{\cdot}_{\infty}$ is  $\norm{\cdot}_1$.
Therefore, the optimal value of this optimization problem is simply $\gamma \norm{h_i}_1 \cdot \norm{x}_p^d$.
%

Now consider the optimization problem:
\begin{equation}\label{eq:nonlinear socp}
\begin{aligned}
\min_{x,\mu,\eta^+,\eta^-} \quad & c^T x\\
\textrm{s.t.} \quad & h_i^T x \leq b_i \quad i = 1, \dots ,r \\
& \sum_j \mu_j v_j + \sum_{kl} \eta^+_{kl} (\gamma\norm{x_k}_p^d + (y_k)_l )   \\
&\qquad + \sum_{kl} \eta^-_{kl} (\gamma\norm{x_k}_p^d - (y_k)_l )  + \gamma \norm{h_i}_1 \cdot \norm{x}_p^d \leq b_i \quad i = 1, \dots, r \\
& \mu \geq 0, \quad \eta^+ \geq 0, \quad \eta^- \geq 0.
\end{aligned}
\end{equation}
If $d=0$, or if $d=1$ and $p\in\{1,+\infty\}$, then \eqref{eq:nonlinear socp} is a linear program.
Otherwise, the rationality of $p$ ensures that $\norm{x}^d_p$ is \emph{second-order cone representable} \citetext{see \citealp[Sect.~2.3]{bental_nemirovski}; \citealp[Sect.~2.5]{LOBO1998193}}.
This means that \eqref{eq:nonlinear socp} is indeed a second-order cone program.


Let $F \subset \R^n$ denote 
the projection of the feasible set of \eqref{eq:nonlinear socp} onto $x$-space.
We claim that the feasible set of \eqref{eq:nonlinear} equals $F \cup \{ x_1, \dots, x_k\}$.
Indeed, since the vectors $x_k$ are one-step safe measurements, 
we have that $x_k \in S$ and $y_k \in S$.
This implies that $x_k$ is feasible to \eqref{eq:nonlinear}.
Furthermore, for $x \in F \setminus \{ x_1, \dots, x_k \}$,
we have shown that $x$ satisfies the constraints of \eqref{eq:nonlinear bilevel}
if and only if $x$ satisfies the constraints of \eqref{eq:nonlinear socp}.

Therefore, optimizing an objective function over the feasible set of \eqref{eq:nonlinear}
is equivalent to optimizing the same objective function over $F \cup \{x_1, \dots, x_k\}$.
%
}

\subsection{Numerical Example}\label{sec:nonlinear example}
\begin{figure}
\figureconts
{fig:nonlinear}
{\caption{One-step safe learning of a nonlinear system associated with the example in \sectionref{sec:nonlinear example}.}}
{%
\subfigure[Dependence of $S^1_0$ on $\gamma$.][c]{%
\label{fig:nonlinear gamma}
\includegraphics[width=.5\textwidth -.5em]{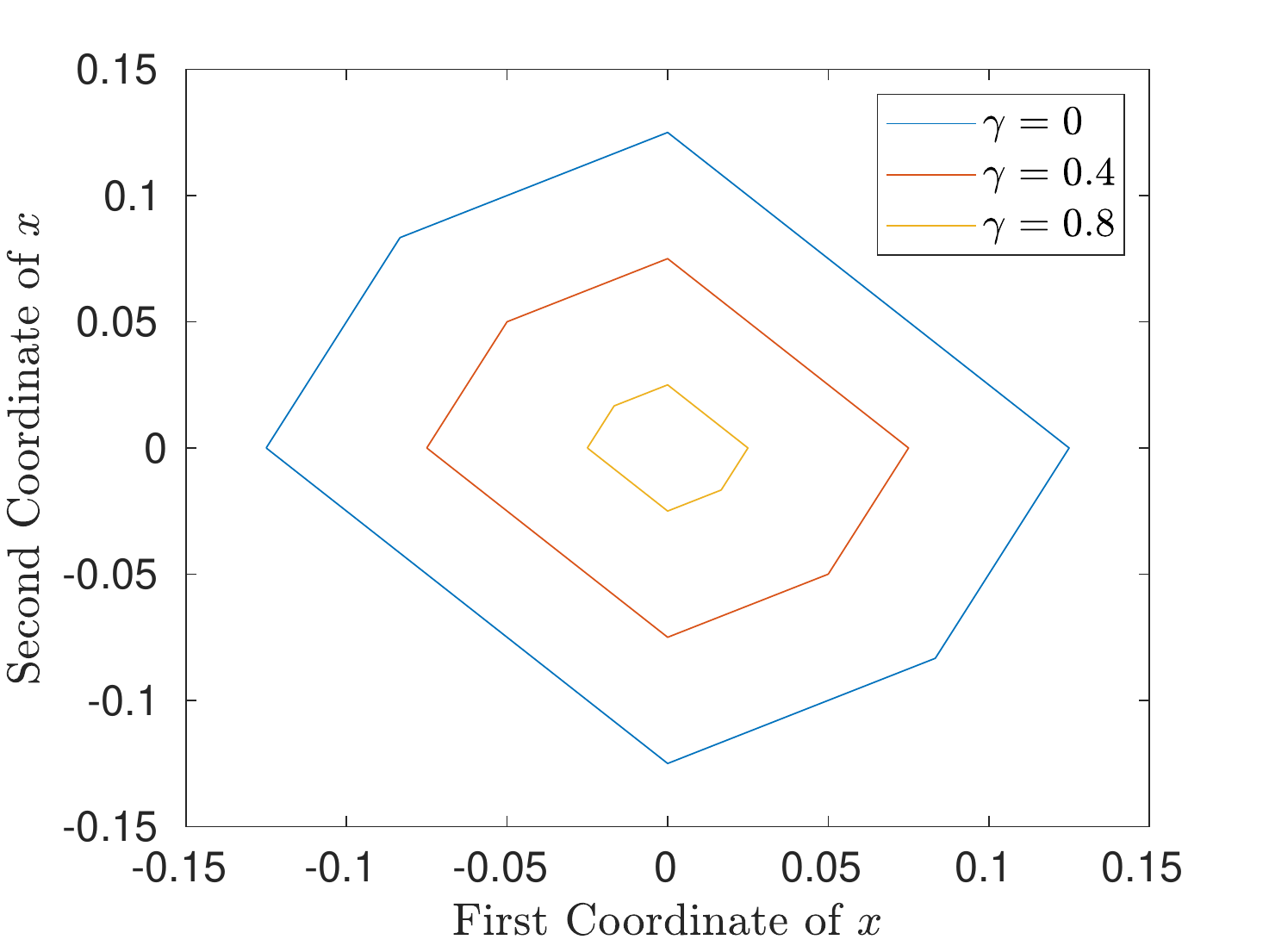}
} 
\subfigure[$S^1_k$ grows with $k$.][c]{%
\label{fig:nonlinear learning}
\includegraphics[width=.5\textwidth -.5em]{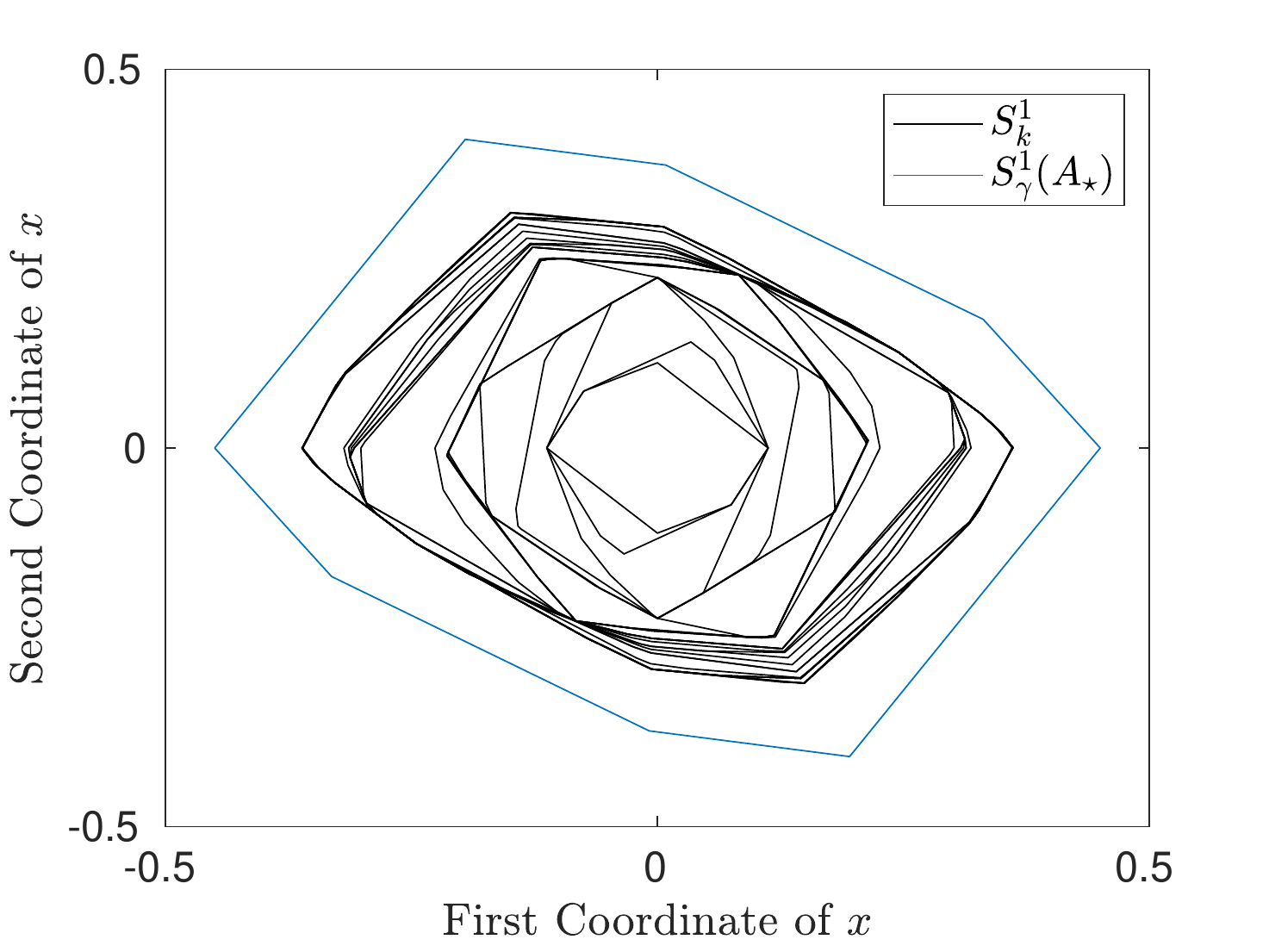}
}
}
\end{figure}

We present a numerical example, again with $n=4$.
Here we take: 
\begin{align*}
      S &= \{ x \in \R^4 \mid |x_i| \leq 1, i=1,\dots,4 \}, \\
    U_{0,A} &= \{ A \in \R^{4 \times 4} \mid -4 \leq A_{ij} \leq 8, i = 1, \dots, 4, j = 1, \dots, 4 \}, \\
    U_{0,g} &= \{ g : \R^4 \rightarrow \R^4 \mid \norm{g(x)}_\infty \leq \gamma \quad \forall x \in S \}.
\end{align*}
In \figureref{fig:nonlinear gamma}, we plot $S^1_0$ (the one-step safety region without any measurements) projected onto the first two dimensions of $x$ for $\gamma \in \{ 0,0.4,0.8 \}$.
As expected, larger values of $\gamma$ result in smaller one-step safety regions.

For our next experiment,
we choose the matrix $\trueA$ in \eqref{eq:nonlinear_model} to be the same matrix used in the example in \sectionref{sec:one-step example}.
We let $\gamma = 0.1$, and
\[ \trueg (x) = \frac{\gamma}{2} 
\left( x_2^2 - x_3 x_4, \quad
\sqrt{x_1^4 + x_3^4}, \quad
x_3 \sin^2(x_1), \quad
\sin^2(x_2) \right) ^T \in U_{0,g}. \]

Since the true system is not linear, we cannot hope to learn the exactly dynamics in $n$ steps
as we did in the linear case.
We instead pick $30$ one-step safe points $x_1,\dots,x_{30}$ (by sequentially solving the second-order cone program from \theoremref{thm:nonlinear socp}) and observe $y_k = \truef(x_k)$ for each $k=1,\dots,30$.
In order to encourage exploration of the state space, we optimize in random directions in every iteration (instead of optimizing the same cost function throughout the process).
In \figureref{fig:nonlinear learning}, we plot $S^1_k$ (the one-step safety region after $k$ measurements) projected onto the first two dimensions of $x$ for $k=0,\dots,30$.
We also plot the projection of $S^1_\gamma (\trueA)$, which we define
as the set of one-step safe points if we knew $\trueA$, but not $\trueg$:
\[ S^1_\gamma (\trueA) \defn \{x \in S \mid \trueA x + g(x) \in S \quad \forall g \in U_{0,g} \}.\]
Finally, we undertake the task of learning the unknown nonlinear dynamics.
We only use information from our first $8$ data points in order to make the fitting task more challenging.
We fit a function of the form
\[\hat{f}(x) = \hat{A}x + \hat{g}(x) , \]
where $\hat{A} \in \R^{4 \times 4}$ and each entry of $\hat{g} : \R^4 \rightarrow \R^4$ 
is a homogeneous quadratic function of $x$.
Our regression is done by minimizing the least-squares loss function
\[ L(\hat{f}) =  \sum_{k=1}^8 \norm{\hat{f}(x_k) - y_k}^2. \]
We train two models.
The first model, $\hat{f}_{\mathrm{ls}}$, minimizes the least-squares loss with no constraints.
The second model, $\hat{f}_{\mathrm{SOS}}$, minimizes the least-squares loss
subject to the constraints that $\hat{A} \in U_{0,A}$, $\norm{\hat{A}x_k - y_k}_{\infty} \leq \gamma$ for $k=1,\dots,8$, and $\hat{g} \in U_{0,g}$.
%
%
The constraint that $\hat{g} \in U_{0,g}$ is imposed via sum of squares constraints \citetext{see, e.g., \citealp{pablothesis,pmlr-v120-ahmadi20a} for details}.
More specifically, we require that for $j=1, \dots, 4$,
\[
\gamma \pm \hat{g}_j(x) = \sigma_0^{j,\pm}(x) + \sum_{i=1}^r \sigma_i^{j,\pm}(x) (b_i - h_i^T x) \quad \forall x \in \R^4.
\]
Here, $\hat{g}_j(x)$ is the $j$-th entry of the vector $\hat{g}(x)$, and the functions $\sigma_i^{j,\pm}$, for $i=0, \dots, r$ and $j=1,\dots,4$, are sum of squares quadratic functions of $x$.
These constraints can be imposed by semidefinite programming.

%
We sample test points $z_1,\dots,z_{1000}$ uniformly at random in $S$ in order to estimate the 
generalization error.
The root-mean-square error (RMSE) is computed as:
\[ \text{RMSE} (\hat{f}) = \sqrt{ \frac{1}{1000} \sum_{j=1}^{1000} \norm{ \hat{f}(z_i) - \truef(z_i)}^2 }. \]
The $\text{RMSE}(\hat{f}_{\mathrm{SOS}})$ of the constrained model is $0.0851$ and the $\text{RMSE}(\hat{f}_{\mathrm{ls}})$ of the unconstrained model is $0.2567$.
We see that imposing prior knowledge with sum of squares constraints
results in a significantly better fit.

\acks{AAA and AC were partially supported by the MURI award of the AFOSR, the DARPA Young Faculty Award, the
CAREER Award of the NSF, the Google Faculty Award, the Innovation Award of the School of Engineering and Applied Sciences at Princeton University, and the Sloan Fellowship.}

\bibliography{paper}

\end{document}